\documentclass{article}

\title{\LARGE \textbf{A Note on Hamilton Cycles}}
\author{Zh.G. Nikoghosyan\footnote{G.G. Nicoghossian (up to 1997)}\\ \\
Institute for Informatics and Automation Problems\\ National Academy of Sciences\\
P. Sevak 1, Yerevan 0014, Armenia\\ E-mail: zhora@ipia.sci.am}

\begin{document}

\maketitle

\begin{abstract}
If $G$ is a more than one tough graph on $n$ vertices with $\delta\ge \frac{n}{2}-a$ for a given $a>0$ and $n$ is large enough then $G$ is hamiltonian. \\


\end{abstract}

\section{Introduction}

Only finite undirected graphs without loops or multiple edges are considered. We use $n$, $\delta$, $c$ and $\tau$ to denote the number of vertices (order),  the minimum degree, circumference and the toughness of a graph, respectively. Let $s(G)$ denote the number of components of a graph $G$. A graph $G$ is $t$-tough if $|S|\ge ts(G\backslash S)$ for every subset $S$ of the vertex set $V(G)$ with $s(G\backslash S)>1$. The toughness of $G$, denoted $\tau(G)$, is the maximum value of $t$ for which $G$ is $t$-tough (taking $\tau(K_n)=\infty$ for all $n\ge 1$).  A good reference for any undefined terms is \cite{[2]}. 

The earliest degree condition for a graph to be hamiltonian is due to Dirac \cite{[3]}.\\

\noindent\textbf{Theorem A \cite{[3]}}. Every graph with $\delta\ge \frac{1}{2}n$ is hamiltonian.\\

Jung \cite{[4]} proved that the minimum degree bound $\frac{1}{2}n$ in Theorem A can be slightly relaxed when $n\ge 11$ and $\tau \ge 1$.\\

\noindent\textbf{Theorem B \cite{[4]}}. Every graph with $n\ge11$, $\tau\ge1$ and $\delta\ge \frac{n}{2}-2$ is hamiltonian.\\

This bound $\frac{n}{2}-2$ itself was lowered to $\frac{n}{2}-3.5$ under stronger conditions $n\ge30$ and $\tau>1$.\\

\noindent\textbf{Theorem C \cite{[1]}}. Every graph with $n\ge30$, $\tau>1$ and $\delta\ge \frac{n}{2}-3.5$ is hamiltonian.\\

In this note we show the following.\\

\noindent\textbf{Theorem 1}. Let $G$ be a graph with $\tau>1$ and $\delta\ge\frac{n}{2}-a$ for a given $a>0$. Then $G$ is hamiltonian if $n$ is large enough.\\

The proof of Theorem 1 easily follows from the following theorem due to Jung and Wittmann \cite{[5]}.\\

\noindent\textbf{Theorem D \cite{[5]}}. In every 2-connected graph, $c\ge\min\{n,(\tau+1)\delta+\tau\}$.\\

\noindent\textbf{Proof of Theorem 1}. Let $\tau=1+\epsilon$ for some $\epsilon>0$. If
$$
n\ge\frac{4a-2}{\epsilon}+2a-2   \eqno{(1)}
$$
then equivalently
$$
(2+\epsilon)\left(\frac{n}{2}-a+1\right)-1\ge n
$$
and by Theorem D,
$$
c\ge(\tau+1)\delta+\tau\ge(2+\epsilon)\left(\frac{n}{2}-a+1\right)-1\ge n.
$$
So, if (1) then $G$ is hamiltonian. Theorem 1 is proved.

\end{document}